\documentclass[12pt]{article}
\input epsf      
\usepackage{graphicx}
\usepackage{amssymb}
\usepackage{amsfonts}
\usepackage{amsthm}
\usepackage{amsfonts}

\newcommand{\complex}{\mathbb{C}}

\newcommand{\half}{\mathbb{H}}

\newcommand{\compose}{\circ}

\newcommand{\blength}{b}
\newcommand{\nterms}{n}

\newcommand{\ba}[1]{\begin{array}{#1}}
\newcommand{\ea}{\end{array}}

\newcommand{\be}{\begin{equation}}
\newcommand{\ee}{\end{equation}}
\newcommand{\bea}{\begin{eqnarray}}
\newcommand{\eea}{\end{eqnarray}}
\newcommand{\beann}{\begin{eqnarray*}}
\newcommand{\eeann}{\end{eqnarray*}}

\def\reff#1{(\ref{#1})}
\begin{document}

\bibstyle{ams}

\title{A Fast Algorithm for Simulating 
\\ the Chordal Schramm-Loewner Evolution}

\author{Tom Kennedy \\
Department of Mathematics \\
University of Arizona \\
Tucson, AZ 85721 \\
http://www.math.arizona.edu/$\, \hbox{}_{\widetilde{}}$ tgk \\
email: tgk@math.arizona.edu \\
}

\maketitle

\begin{abstract}
The Schramm-Loewner evolution (SLE) can be simulated by dividing the
time interval into $N$ subintervals and approximating the random conformal
map of the SLE by the composition of $N$ random, but relatively simple,
conformal maps. In the usual implementation the time required to compute
a single point on the SLE curve is $O(N)$. We give an algorithm for 
which the time to compute a single point is $O(N^p)$ with $p<1$. 
Simulations with $\kappa=8/3$ and $\kappa=6$ both give a value of 
$p$ of approximately $0.4$. 
\end{abstract}

\newpage

\section{Introduction}

The Schramm-Loewner evolution (SLE) is a stochastic process that produces 
a random curve in the complex plane. In this paper we will be concerned 
with chordal SLE in which the random curve, the SLE ``trace,'' lies in 
the upper half plane and goes from $0$ to $\infty$.  
The classical Loewner equation is 
\be 
\partial_t \, g_t(z) = {2 \over g_t(z) - U_t}  
\label{loweq} 
\ee
with the initial condition $g_0(z)=z$. 
Here $U_t$ is a real-valued ``driving function.'' The equation 
defines a one parameter family $g_t$ of conformal maps from a simply 
connected subset of the upper half plane onto the upper half plane. 
SLE is obtained by taking $U_t = \sqrt{\kappa} B_t$ where $\kappa>0$
is a parameter and $B_t$ is a Brownian motion with mean zero and 
variance $t$. For further discussion of SLE we refer the reader to 
\cite{lawler,werner} and the original references \cite{rs,schramm}.

The most common method for simulating SLE is not to numerically solve
the above differential equation. Instead one partitions the 
time interval into $N$ subintervals and approximates
$g_t$ by the composition of a sequence of $N$ conformal maps which are 
approximations to the solution of the Loewner equation over the 
subintervals. Computing a point on the SLE trace requires evaluating 
the composition of roughly $N$ conformal maps and so takes a time $O(N)$. 
We will refer to the time it takes to compute one point on the SLE trace
as the ``time per point.'' 

\begin{figure}[htb]
\includegraphics{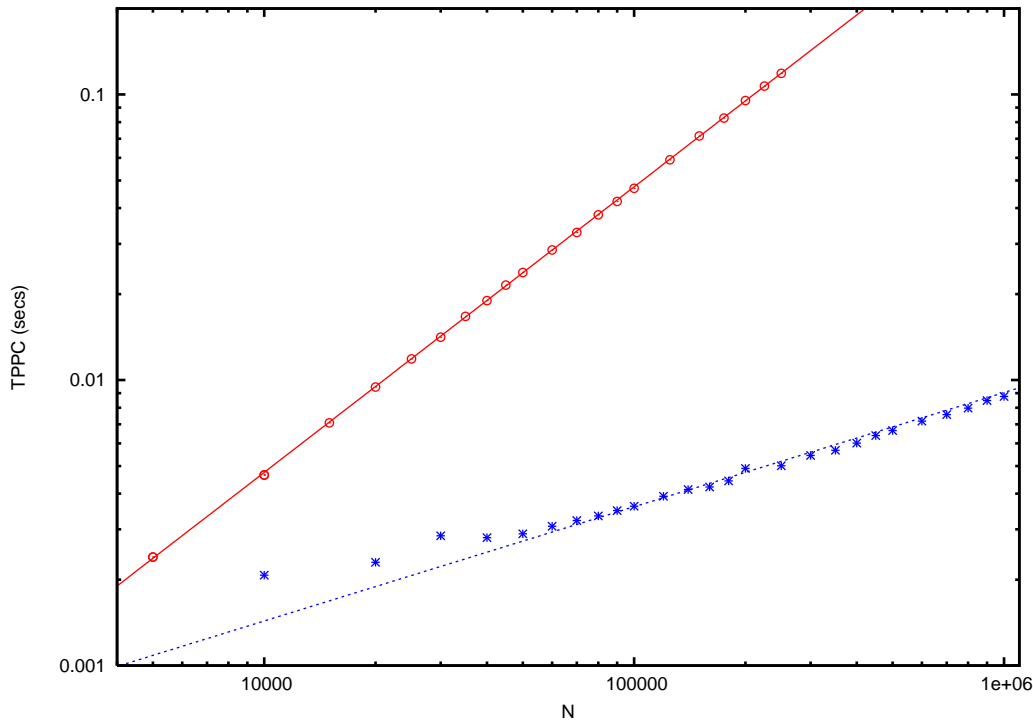}
\caption{ Time per point computed as a function of $N$, the number of
subintervals in the partition of the time interval for $\kappa=8/3$. 
The top curve is the usual algorithm; the bottom curve is the new algorithm.
The lines shown have slopes $1$ and $0.4$. }
\label{tppcprelim_8_3}
\end{figure}

\begin{figure}[htb]
\includegraphics{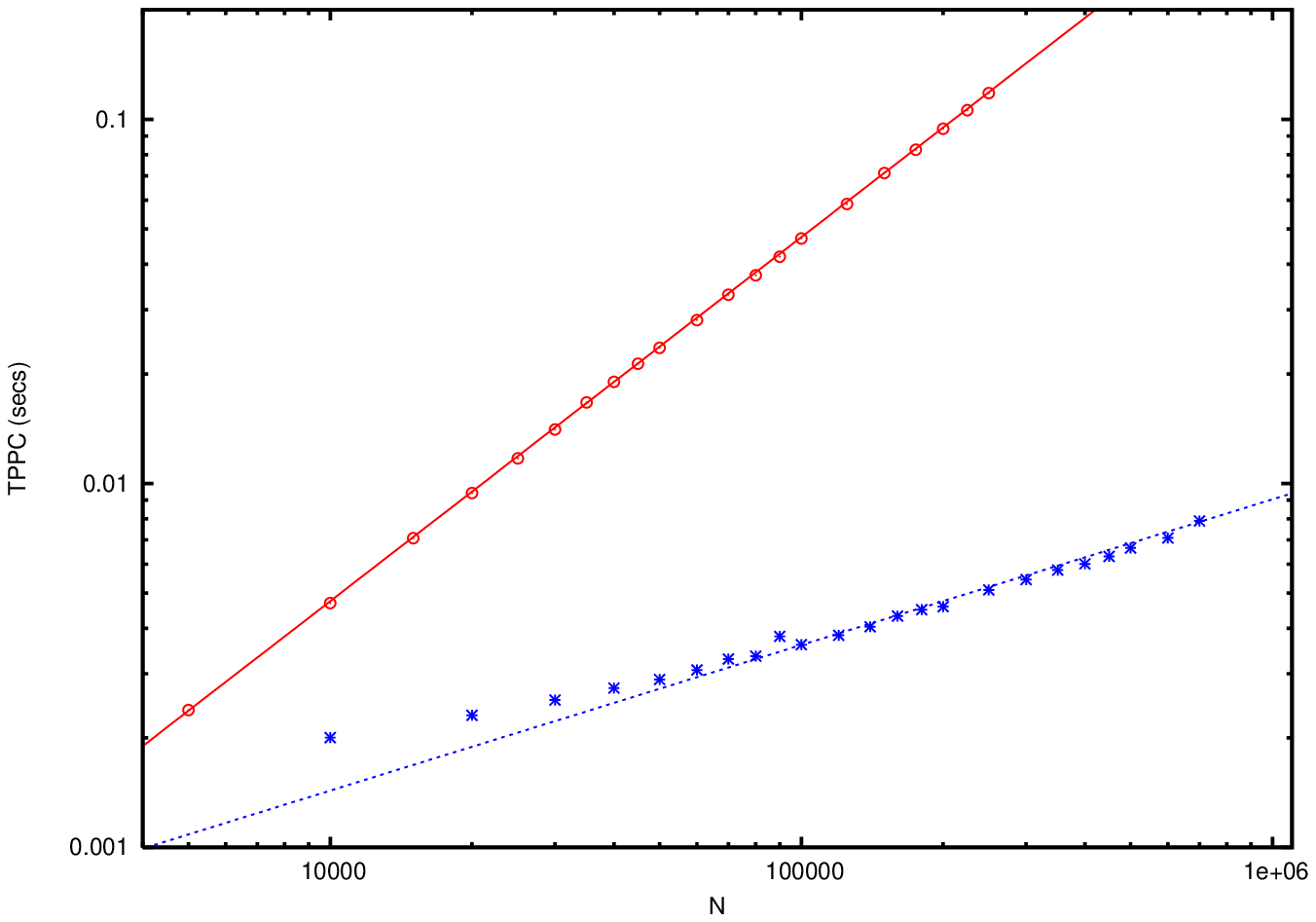}
\caption{ Time per point computed as a function of $N$, the number of
subintervals in the partition of the time interval for $\kappa=6$.  
The lines shown have slopes $1$ and $0.4$. }
\label{tppcprelim_6}
\end{figure}

The goal of this paper is to give an algorithm for which the time per 
point is $O(N^p)$ with $p<1$. We do not prove that our algorithm does this. 
The time our algorithm takes to generate an SLE curve depends on the 
behavior of the particular curve. For certain atypical curves, our 
algorithm will not be faster than the usual algorithm. 
So a rigorous analysis of our algorithm is a daunting task. 
We have studied the behavior of the algorithm by simulation for 
$\kappa=8/3$ and $\kappa=6$.
For both of these values of $\kappa$, the time per point for 
our algorithm is $O(N^p)$ with $p$ around $0.4$. 
Figures \ref{tppcprelim_8_3} and \ref{tppcprelim_6} 
show log-log plots of the time per point as a function of $N$ for 
the usual algorithm and for our new algorithm. In the first figure
$\kappa=8/3$, and in the second it is $6$. Although the 
behavior of the SLE trace is qualitatively different for these two 
values of $\kappa$, the behavior of the time required for the simulation
is quite similar. 
The straight lines show fits by $c_1 N$ for the 
usual algorithm and by $c_2 N^{0.4}$ for our algorithm.
For $N=100,000$ our algorithm is faster by approximately a factor of $14$ 
and for $N=1,000,000$ our algorithm is faster by 
approximately a factor of $56$. 

%

%

There is an inverse problem that is closely related to the present paper.
Given a simple curve in the upper half plane, compute the corresponding
driving function in the Loewner equation.
An important application of this inverse problem is studying if 
a family of random curves is SLE. Given a large sample of the curves, 
one computes the corresponding samples of the driving function and 
then tests if they are a Brownian motion. This was done for 
interfaces in two-dimensional spin glass ground states in \cite{ahhm,bdm}
and for certain isolines in two-dimensional models of turbulence in 
\cite{bbcfA, bbcfB}. 
The methods of this paper apply to this inverse problem as well. The naive
implementation of the zipper algorithm to compute the driving function 
for a curve with $N$ points runs in a time $O(N^2)$.
Using the methods of this paper, the algorithm runs in a time $O(N^{1.35})$.
This fast algorithm is studied in \cite{ken}.

In section \ref{sectdiscretize} we explain the standard method for 
``discretizing'' SLE. A particular example was studied by Bauer 
\cite{bauer}. The material in this section is well known but does not
seem to have appeared in the literature yet.
In section \ref{sectlaurent} we explain our new algorithm. 
The algorithm involves several parameters, so in section 
\ref{sectparms} we study how the error and time required for the 
algorithm depend on these parameters. 
The appendix gives some power series needed in particular implementations of 
the algorithm. C++ code implementing the new algorithm may be 
downloaded from the author's homepage.

\section{Discretizing SLE}
\label{sectdiscretize}

We are going to approximate SLE by a 
``discrete SLE.'' The term discrete SLE is a bit misleading. 
The random process we define produces continuous curves in the 
upper half plane. The plane is not replaced by a lattice. 
Let $0=t_0 < t_1 < t_2 < \cdots < t_n=1$ be a partition of the time
interval $[0,t]$. (Thanks to the scaling property of SLE, it is no loss
of generality to take the time interval to be $[0,1]$.)
The times $t_k$  play a special role, but the random curves 
are still defined for all time. 
We can think of our discrete approximation as the result of replacing
the Brownian motion in the driving function by some stochastic
process that approximates (or even equals) the Brownian motion at the 
times $t_k$, and is defined in between these times so that the Loewner 
equation may be solved explicitly. 

We begin by reviewing some facts about the Schramm-Loewner evolution for the 
upper half-plane,
\be 
\partial_t \, g_t(z) = {2 \over g_t(z) - \sqrt{\kappa} B_t}  
\ee
with $g_0(z)=z$. 
The set $K_t$ contains the points $z$ in the upper half-plane for which the 
solution to this equation no longer exists at time $t$. 
Let $t,s>0$. The map $g_{t+s}$ maps $\half \setminus K_{t+s}$ onto $\half$. 
We can do this in two stages. We first apply the map $g_s$. 
This maps $\half \setminus K_{s}$ onto $\half$, and it maps  
$\half \setminus K_{t+s}$ onto $\half \setminus g_s(K_{t+s} \setminus K_s)$. 
Let $\hat{g}_t$  be the conformal map that maps 
$\half \setminus g_s(K_{t+s} \setminus K_s)$ onto $\half$ with the usual 
hydrodynamic normalization. By the uniqueness of these maps, 
\be
g_{s+t}=\hat{g}_t \compose g_s, \quad i.e., \quad 
\hat{g}_t = g_{s+t} \compose g_s^{-1}
\ee
Then 
\be
{d \over dt} \hat{g}_t(z)
={d \over dt} g_{s+t} \compose g_s^{-1}(z) 
= {2 \over g_{s+t} \compose g_s^{-1}(z) -U_{s+t}}
= {2 \over \hat{g}_{t}(z)-U_{s+t}}
\ee
Note that $\hat{g}_0(z)=z$. Thus  $\hat{g}_t(z)$ is obtained by 
solving the Loewner equation with driving function $\hat{U}_t=U_{s+t}$.
This driving function starts at $U_s$, and so the 
$\hat{K}_t$ associated with $\hat{g}_t$ starts growing at $U_s$.
In our approximation the driving function will be sufficiently nice that
$\hat{K}_t$ is just a curve which starts at $U_s$.  

We now return to the partition $0=t_0 < t_1 < t_2 < \cdots t_n =1$, and 
define
\be
G_k = g_{t_k} \compose g_{t_{k-1}}^{-1}
\ee
So 
\be 
g_{t_k} = G_k \compose G_{k-1} \compose G_{k-2} \compose \cdots \cdots 
G_2 \compose G_1 
\ee
By the remarks above, $G_k(z)$ is obtained by solving the Loewner equation
with driving function $U_{t_{k-1} +t}$ for $t=0$ to $t=\Delta_k$,
where $\Delta_k = t_k-t_{k-1}$. 
Note that $G_k$ maps 
$\half$ minus a set ``centered'' around $U_{t_{k-1}}$ to $\half$.
If we consider $G_k(z+U_{t_{k-1}+t}) - U_{t_{k-1}}$, 
it is obtained by solving the 
Loewner equation with driving function $U_{t_{k-1}+t} - U_{t_{k-1}}$
for $t=0$ to $t=\Delta_k$. This driving function starts at 
$0$ and ends at $\delta_k$ where $\delta_k=U_{t_k}-U_{t_{k-1}}$.
So this conformal map takes 
$\half$ minus a set ``centered'' around the origin onto $\half$.

The approximation to the SLE trace is given by 
$\gamma(t)=g_t^{-1}(U_t)$.
Let $z_k=g_{t_k}^{-1}(U_{t_k})$.
We will only consider the points on this curve which correspond to 
times $t=t_k$. One could consider other points on the curve, 
but the distance between consecutive $z_k$ is already of the order of the 
error in our approximation, so there is little reason to consider more points. 
The points $z_k$ are given by 
\be 
z_k= G^{-1}_1 \compose G^{-1}_2 \compose \cdots \cdots 
G^{-1}_{k-1} \compose G^{-1}_k (U_{t_k}) 
\ee
We would like to rewrite this using conformal maps that 
depend only on the change in $U_t$ over the time intervals. 
Define
\be
h_k(z)=G^{-1}_k(z+U_{t_k})-U_{t_{k-1}}
\ee
So 
\be
z_k= h_1 \compose h_2 \compose \cdots \cdots 
\compose h_{k-1} \compose h_k(0)
\label{dsle}
\ee
Recall that if we solve the Loewner equation with driving function 
$U_{t_{k-1}+t} - U_{t_{k-1}}$ for $t=0$ to $t=\Delta_k$,
we get $g_k(z)$ where 
\be
g_k(z)=G_k(z+U_{t_{k-1}}) - U_{t_{k-1}}
\ee
Letting $f_k(z)=g^{-1}_k(z)$, we have 
\be
f_k(z)=G^{-1}_k(z+U_{t_{k-1}})-U_{t_{k-1}}
\ee
and so $h_k(z)=f_k(z+\delta_k)$. 
As noted above, $f_k$ maps $\half$ to $\half$ minus a curve that starts at $0$.
The driving function ends at $\delta_k$, so $f_k(\delta_k)$ is the tip 
of the curve. 
It follows that $h_k(z)$ maps $\half$ onto $\half$ minus a curve 
starting at $0$ and $h_k(0)$ is the tip of the curve.
Thus we have the following simple picture for eq. \reff{dsle}.
The first map $h_k$ welds together a small interval on the real 
axis containing the origin to produce a small cut. The origin is mapped
to the tip of this cut. The second map welds together 
a (possibly different) small interval in such a way that it produces another
small cut. The original cut is moved away from the origin with its base 
being at the tip of the new cut. This process continues. 
Each map introduces a new small cut whose tip is attached to 
the image of the base of the previous cut. 

The key idea is to define $U_t$ for $t_{k-1} \le t \le t_k$ so that 
$g_k(z)$ may be explicitly computed. 
There are two constraints
on $g_k$. The curve must have capacity $2 \Delta_k$ and 
$g_k$ must map the tip of the curve to $\delta_k$. 
Given any simple curve satisfying these two constraints and starting at the 
origin, it will be the solution of the Loewner equation for some 
driving function which goes from $0$ to $\delta_k$ over the time 
interval $[0,\Delta_k]$. Different choices of this interpolating curve give 
us different discretizations. Before we discuss particular 
discretizations, we will discuss the definitions of $\Delta_k$ and $\delta_k$. 

The simplest choice for $\Delta_k$ is to use a uniform partition of the 
time interval, $\Delta_k=1/N$. However, if we use a uniform partition and 
look at the resulting points $z_k$, they appear to be farther apart 
on average at the beginning of the curve than at the end.
To understand this, consider the case of $\kappa=8/3$ which is believed 
to correspond to the self-avoiding walk. Let $\omega(s)$ denote the 
scaling limit of the self-avoiding walk. 
With its natural parameterization we have
$E \, \omega(s)^2 = c s^{2 \nu}$ where $\nu=3/4$ and $c$ is a constant. 
If we take points at equally spaced times using this parameterization we 
will get points on the walk that are roughly equally distant. 
For $SLE$ with its parameterization using capacity, 
we have $E \, \omega(t)^2 = c t$ for some constant $c$. 
To match these two parameterizations in an average sense we should 
take $t=s^{2 \nu}$. Thus to get points approximately equally spaced 
on the SLE curve, we should take $t_k=(k/N)^{2 \nu}$.

The $\delta_k$ should be chosen so that the stochastic process $U_t$
will converge to $\sqrt{\kappa}$ times Brownian motion as 
$N \rightarrow \infty$. 
One possibility is take the $\delta_k$ to be independent normal random 
variables with mean zero and variance $\kappa \Delta_k$. If we do this, 
then $U_t$ and $\sqrt{\kappa} B_t$ will have the same distributions if 
we only consider times chosen from the $t_k$. 
Another possibility is to take 
the $\delta_k$ to be independent random variables with 
$\delta_k = \pm \sqrt{\kappa \Delta_k}$ where the choices of $+$ and $-$ are 
equally probable. If we use this choice with the uniform partition of the 
time interval, then we are approximating the Brownian motion by a simple 
random walk. 

We now consider specific choices of the interpolating curve used for the 
discretization. A popular choice is the following.
Let $C$ be a line segment starting at the origin with a polar 
angle of $\alpha \pi$. 
$g_k$ maps $\half \setminus C$ onto $\half$. There are two 
degrees of freedom for the line segment - its length and $\alpha$. 
There are two constraints - the line segment must have capacity $2 \Delta_k$
and the tip of the segment must get mapped to $\delta_k$. 

Consider the map 
\be
\phi(z)=(z+y)^{1-\alpha} (z-x)^\alpha
\ee
where $x,y>0$. 
It maps the half plane onto the half plane minus a line segment which 
starts at the origin and forms an angle $\alpha$ with the positive real 
axis. The interval $[-y,x]$ gets mapped onto the slit. The length of this
interval determines the length of the slit. 
Shifting this interval (relative to $0$) does not change the length of the 
slit. To obtain the map $g_k$, we must choose $x$ and $y$ so that $g_k$
satisfies the hydrodynamic normalization and has capacity $2 \Delta_k$. 
Tedious but straightforward calculation shows if we let 
\be
f_t(z)=
\left(z+ 2 \sqrt{t} \sqrt{\alpha \over 1-\alpha} \right)^{1-\alpha} 
\left(z- 2 \sqrt{t} \sqrt{1-\alpha \over \alpha} \right)^\alpha
\label{tiltf}
\ee
then $f_t^{-1}(z)$ satisfies the hydrodynamic normalization 
and has capacity $2 t$. In particular, $f_k(z)$ is given by the above 
equation with $t=\Delta_k$. 
We know from the general theory that $g_t(z)=f_t^{-1}(z)$ 
satisfies the Loewner equation \reff{loweq} 
for some driving function $U_t$. Some calculation then shows that 
\be
U_t = c_\alpha \sqrt{t}
\ee
where 
\be
c_\alpha=2 {1-2 \alpha \over \sqrt{\alpha (1-\alpha)}}
\ee
For the map $g_k$ we need $U_{t_k}-U_{t_{k-1}}=\delta_k$, and so 
\be
\delta_k = c_\alpha \sqrt{\Delta_k}
\label{alphaeq} 
\ee
This equation determines $\alpha$. ($\alpha$ depends on $k$.)

Define 
\be
v={\delta_k^2 \over \Delta_k}
\ee
Then squaring \reff{alphaeq} gives
\be
c_\alpha^2 = v
\ee
which leads to 
\be
16 \alpha^2 + v \alpha^2 - 16 \alpha - v \alpha +4 =0
\ee
and so 
\be
\alpha = {1 \over 2} \pm {1 \over 2} \sqrt{v \over 16 + v}
\ee
We take the choice with $\alpha<1/2$ if $\delta_k>0$ and the choice 
with $\alpha>1/2$ if $\delta_k<0$.
Using $h_k(z)=f_k(z+\delta_k)$ we find that 
\be
h_k(z)=
\left(z+ 2 \sqrt{\Delta_k (1-\alpha) \over \alpha} \right)^{1-\alpha} 
\left(z- 2 \sqrt{\Delta_k \alpha \over 1-\alpha} \right)^\alpha
\ee
Note how confusingly similar this formula is to \reff{tiltf}.

\bigskip

Another discretization is to take 
\be
h_k(z)=\sqrt{z^2-4 \Delta_k} + \delta_k
\ee
This conformal map produces a vertical slit based at $\delta_k$ with 
capacity $2 \Delta_k$. 
It does not map the origin to the tip of the slit. So composing these
maps does not produce a curve. Nonetheless, this discretization converges 
to SLE \cite{bauer}. A vertical slit corresponds to a constant driving 
function. So this discretization corresponds to replacing the Brownian 
motion by a stochastic process that is constant on each time subinterval 
and jumps discontinuously at the times $t_k$. 

\section{A faster algorithm}
\label{sectlaurent}

To motivate what we do in this section, we first consider the speed 
of the algorithm described in the previous section. Recall that 
points on the approximation to the SLE trace are given by 
\be
z_k= h_1 \compose h_2 \compose \cdots \cdots 
\compose h_{k-1} \compose h_k(0)
\label{compose}
\ee
The number of operations needed to compute a single $z_k$ is
proportional to $k$. So to compute all the points $z_k$ with $k=1,2,\cdots N$
requires a time $O(N^2)$. 

It is important to note that the computation of $z_k$ does not depend on 
any of the other $z_j$. So we can compute a subset of the points $z_k$
if we desire. (As an extreme example, if we are only interested in 
$z_N=\gamma(1)$, the time required for the computation is $O(N)$ not
$O(N^2)$.) For the timing tests in this paper we compute the points $z_{jd}$ 
with $j=1,2,\cdots,N/d$ where $d$ is some integer. 
But we emphasize that our algorithm works for any choice of the set of 
points to compute. 
For the above algorithm the time grows as $N^2/d$. The time
per point grows as $N$. We use the time per point throughout
this paper to study the efficiency. It is a natural measure since 
it depends on how finely we discretize the time interval but not on 
the number of points we choose to compute.
The total time to compute the SLE trace is given by the number of 
points we want to compute on it times the time per point.
Our goal is to develop an algorithm for which the time per point 
is $O(N^p)$ with $p<1$. 

Our algorithm begins by grouping the functions in \reff{compose} into
blocks. The number of functions in a block will be denoted by $\blength$. 
Let
\be
H_j = h_{(j-1)b+1} \compose h_{(j-1)b+2} \compose \cdots \compose h_{jb}
\label{blockdef}
\ee
If we write $k$ as $k=mb+l$ with $0 \le l < b$, then we have 
\be
z_k = H_1 \compose H_2 \compose \cdots \compose H_m \compose 
h_{mb+1} \compose h_{mb+2} \compose \cdots \compose h_{mb+l}(0)
\label{blockcompose}
\ee
The number of compositions in \reff{blockcompose} is smaller than the 
number in \reff{compose} by roughly a factor of $b$. Unfortunately, 
even though the $h_i$ are relatively simple, the $H_j$ cannot be explicitly 
computed. Our strategy is to approximate the $h_i$ by functions 
whose compositions can be explicitly computed to give an explicit
approximation to $H_j$. For large $z$, $h_i(z)$ is well approximated 
by its Laurent series about $\infty$. One could approximate $h_i$ 
by truncating this Laurent series. This is the spirit of our approach, 
but our approximation is slightly different.

Let $f(z)$ be a conformal map from $\half$ onto $\half \setminus \gamma[0,t]$,
where $\gamma : [0,t] \rightarrow \half$ is a curve in the upper half
plane with $\gamma(0)=0$. 
We assume that $f(\infty)=\infty$, $f^\prime(\infty)=\infty$ and 
$f(0)=\gamma(t)$.
Let $a,b>0$ be such that $[-a,b]$ is mapped onto the slit $\gamma[0,t]$. 
So $f$ is real valued on $(-\infty,-a]$ and 
$[b,\infty)$. 
By the Schwartz reflection principle,
$f$ has an analytic continuation to $\complex \setminus [-a,b]$,
which we will simply denote by $f$.   
Let $R=\max \{ a,b\}$, so $f$ is analytic on $\{z: |z|>R\}$ and maps
$\infty$ to itself. Thus $f(1/z)$ is analytic on $\{z: 0<|z|<1/R\}$ and 
our assumptions on $f$ imply it has a simple pole at the origin with 
residue 1. So we have 
\be
f(1/z)= 1/z + \sum_{k=0}^\infty \, c_k \, z^{k}
\ee
This gives the Laurent series of $f$ about $\infty$. 
\be
f(z)= z + \sum_{k=0}^\infty \, c_k \, z^{-k}
\ee
If we truncate this Laurent series, it will be a good approximation
to $f$ for large $z$. At first sight, this Laurent series is the natural 
approximation to use for $f$. However, we will use a different but 
closely related representation.

Define $\hat{f}(z)=1/f(1/z)$. Since $f(z)$ does not vanish on 
$\{ |z|>R\}$, $\hat{f}(z)$ is analytic in $\{ z : |z| < 1/R \}$.
Our assumptions on $f$ imply that $\hat{f}(0)=0$ and $\hat{f}^\prime(0)=1$.
So $\hat{f}$ has a power series of the form 
\be
\hat{f}(z) = \sum_{j=0}^\infty \, a_j z^j
\label{hps}
\ee
with $a_0=0$ and $a_1=1$. 
It is not hard to show that $1/R$ is the radius of convergence of 
this power series. 
We will refer to this power series as the ``hat power series'' of $f$. 
Note that the coefficients of the hat power series of $f$ are the
coefficients of the Laurent series of $1/f$. 

The primary advantage of this hat power series over the 
Laurent series is its behavior with respect to composition. 
\be
(f \compose g) \, \hat{} \,(z) = 1/f(1/\hat{g}(z)) = \hat{f}(\hat{g}(z))
\ee
Thus 
\be
(f \compose g) \, \hat{} \, = \hat{f} \compose \hat{g}
\label{composeprop}
\ee
Our approximation for $f(z)$ is to replace $\hat{f}(z)$ by the truncation of 
its power series at order $\nterms$. So 
\be
f(z) = {1 \over \hat{f}(1/z)}
\approx \left[ \sum_{j=0}^n \, a_j z^{-j} \right]^{-1}
\ee

For each $h_i$ we compute the power series of $\hat{h_i}$ to order $\nterms$. 
We then use them and \reff{composeprop} to compute the 
power series of $\hat{H_j}$ to 
order $\nterms$. Let $1/R_j$ be the radius of convergence for the
power series of $\hat{H_j}$. ($R_j$ is easy to compute. It is the smallest 
positive number such that $H_j(R_j)$ and $H_j(-R_j)$ are both real.)
Now consider equation \reff{blockcompose}. 
If $z$ is large compared to $R_j$, then $H_j(z)$ is well approximated 
using its hat power series. We introduce a parameter $L>1$ and 
use the hat power series to compute $H_j(z)$ whenever 
$|z| \ge L R_j$. When $|z| < L R_j$, we just use \reff{blockdef} to 
compute $H_j(z)$. The argument of $H_j$ is the result of applying the 
previous conformal maps to $0$, and so is random. Thus whether 
or not we can approximate a particular $H_j$ using its hat power 
series depends on the randomness and on which $z_k$ we are computing.

\section{Choosing the parameters}
\label{sectparms}

Our algorithm depends on three parameters. The integer $\blength$ is 
the number of functions in a block. The integer $\nterms$ is the order 
at which we truncate the hat power series. The real number $L>1$ 
determines when we use the hat power series approximation to the
block function. These three parameters control how good our approximation
is. We can compute the error that arises from using the hat power
series by computing the discretized SLE curve both using the hat 
power series and not using them. We then define the error to be 
the average distance between the points on the two curves. 
Given some desired level of error, we want to minimize 
the time subject to the constraint that the error is within the 
desired tolerance. 
There is a completely different kind of error - that introduced 
by approximating the Brownian motion by some other stochastic process. 
It should converge to zero as $N \rightarrow \infty$.
The nature of this convergence and in particular its
dependence on the method of discretization is an interesting question,
but we do not study it in this paper. 

The behavior of our algorithm is random in that it depends on the behavior 
of the particular SLE sample. Since the behavior of SLE depends 
qualitatively on $\kappa$, one might expect that the behavior of our
algorithm will depend significantly on $\kappa$. We have studied the 
algorithm for $\kappa=8/3$ and $\kappa=6$ and have found that 
the behaviors for these two values of $\kappa$ are remarkably similar.
We will restrict our discussion and our plots to the case of $\kappa=8/3$,
and discuss how $\kappa=6$ compares  at the end of this section. 

We continue to use the time per point (total time divided by the number of 
points computed) as our measure of the speed of the algorithm.
For this new algorithm it is essentially independent of $d$, the 
number of time intervals between consecutive points computed, provided 
$d$ is not huge. The computation of the hat power series 
for the conformal maps $H_j$ does not depend on how many points we compute. 
When $d$ is large enough, the time required for this computation 
will dominate, and the time per point will no longer be independent of $d$. 

\begin{figure}[htb]
\includegraphics{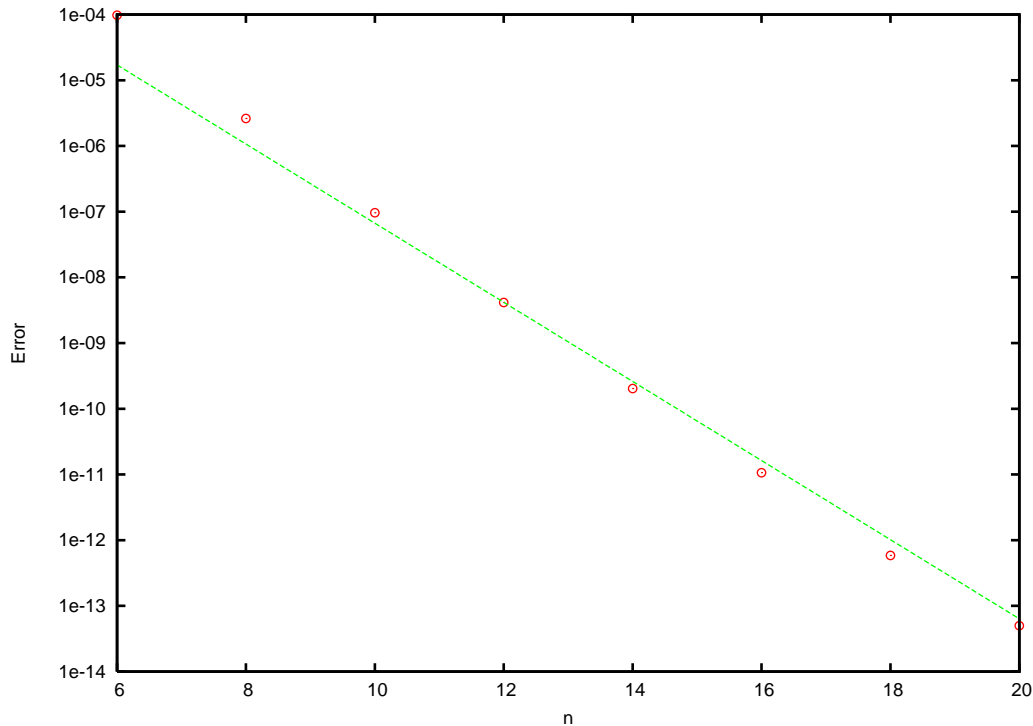}
\caption{ Error as a function of $\nterms$, the order of the hat 
power series. The line is a fit by $c L^{-n}$. 
}
\label{err_nterms}
\end{figure}

\begin{figure}[htb]
\includegraphics{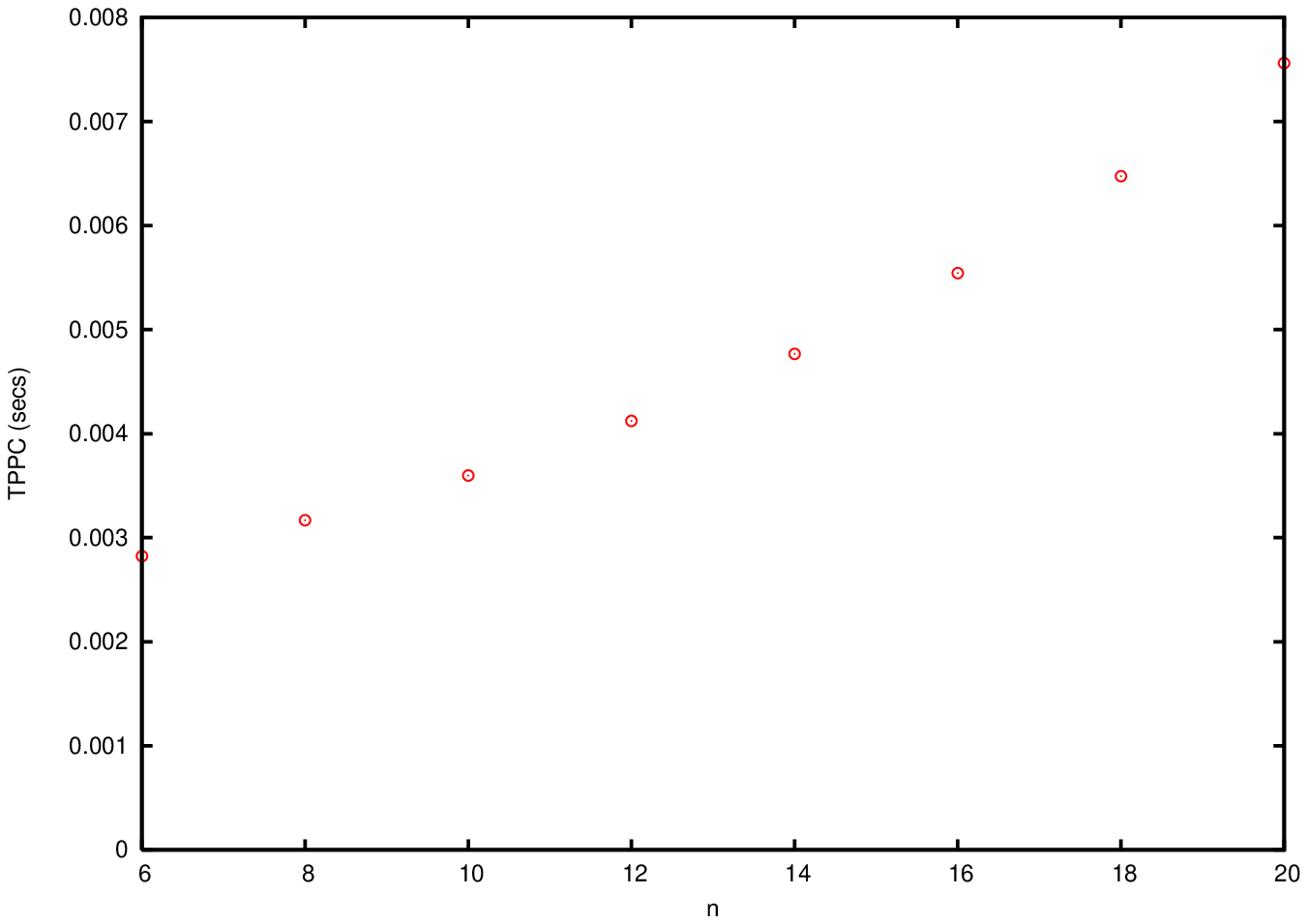}
\caption{ Time per point computed as a function of $\nterms$}
\label{time_nterms}
\end{figure}

We first consider the effect of $\nterms$.
Obviously, as $\nterms$ grows the error should decrease 
and the time should increase.  
Figure \ref{err_nterms} shows that as a function of $\nterms$,
the error is approximately proportional to $L^{-\nterms}$. 
Figure \ref{time_nterms} shows the growth of the 
time with respect to $\nterms$. 
(In both figures we have taken $N=100,000$, $\blength=40$, and $L=4$.)

\begin{figure}[htb]
\includegraphics{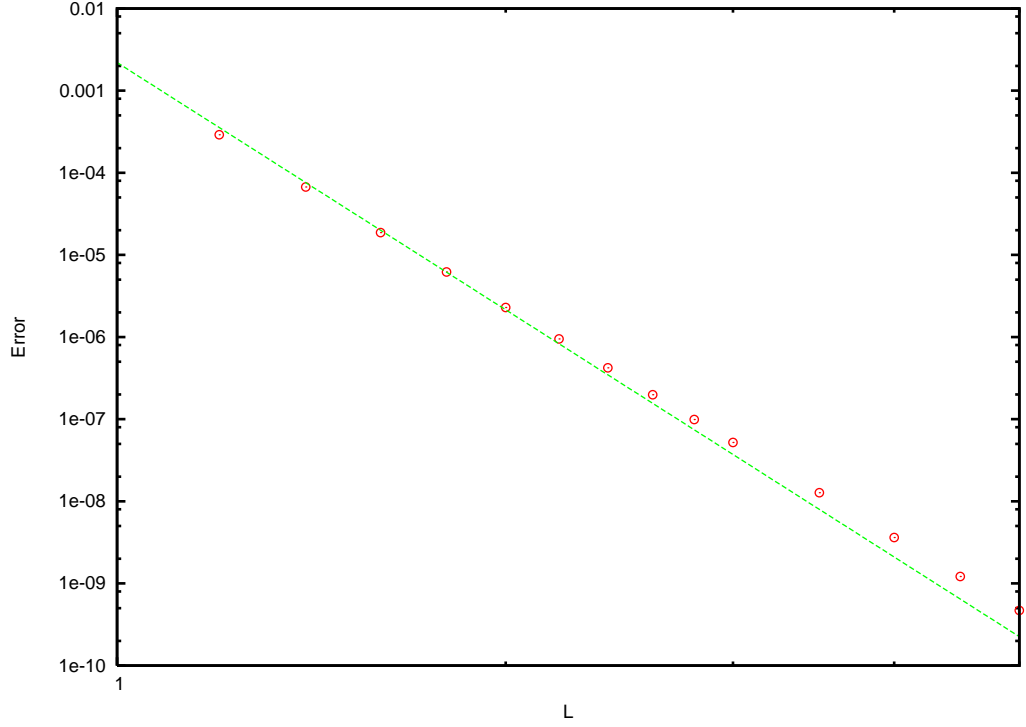}
\caption{ Error as a function of $L$, the parameter that determines how 
often we use the hat power series approximation. 
The line is a fit by $c L^{-n}$. 
}
\label{err_L}
\end{figure}

\begin{figure}[htb]
\includegraphics{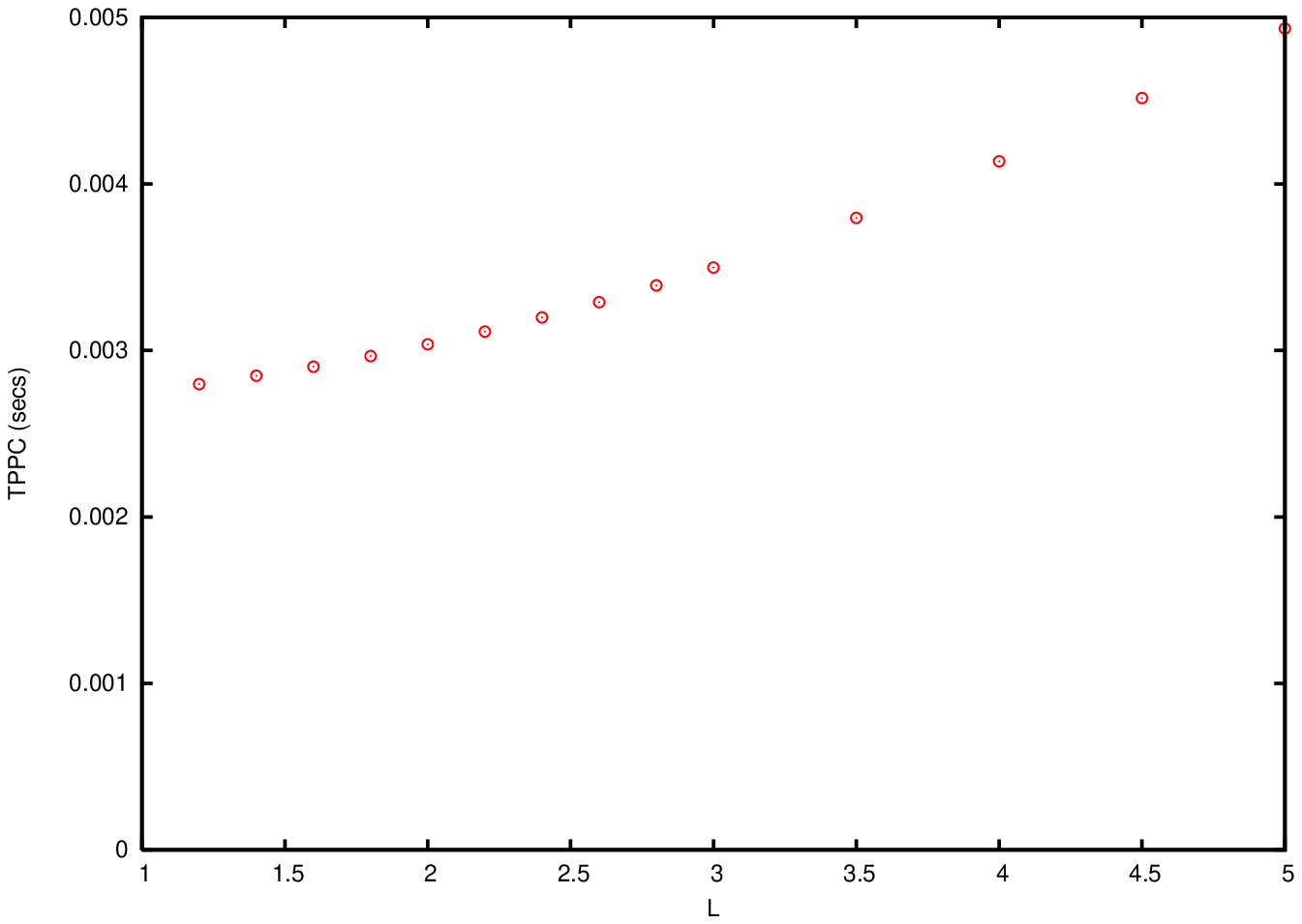}
\caption{ Time per point computed as a function of $L$.}
\label{time_L}
\end{figure}

Next we consider the effect of $L$. As $L$ increases we use the 
hat power series approximation only for larger $z$ and so the 
error should decrease. However, using the approximation less frequently
will increase the time required. 
The error decreases with $L$ roughly as $L^{-\nterms}$. 
(Figure \ref{err_L})
The dependence of the time on $L$ shown in figure \ref{time_L} 
is somewhat surprising. The time does not grow with $L$ as quickly as one 
might expect. Eventually, 
as $L$ gets large the time will be on the order of the time for the 
standard algorithm that does not use hat power series, but 
this only happens at values of $L$ considerable larger than those shown.
(In both figures we have taken $N=100,000$, $\blength=40$, and $\nterms=12$.)

\begin{figure}[htb]
\includegraphics{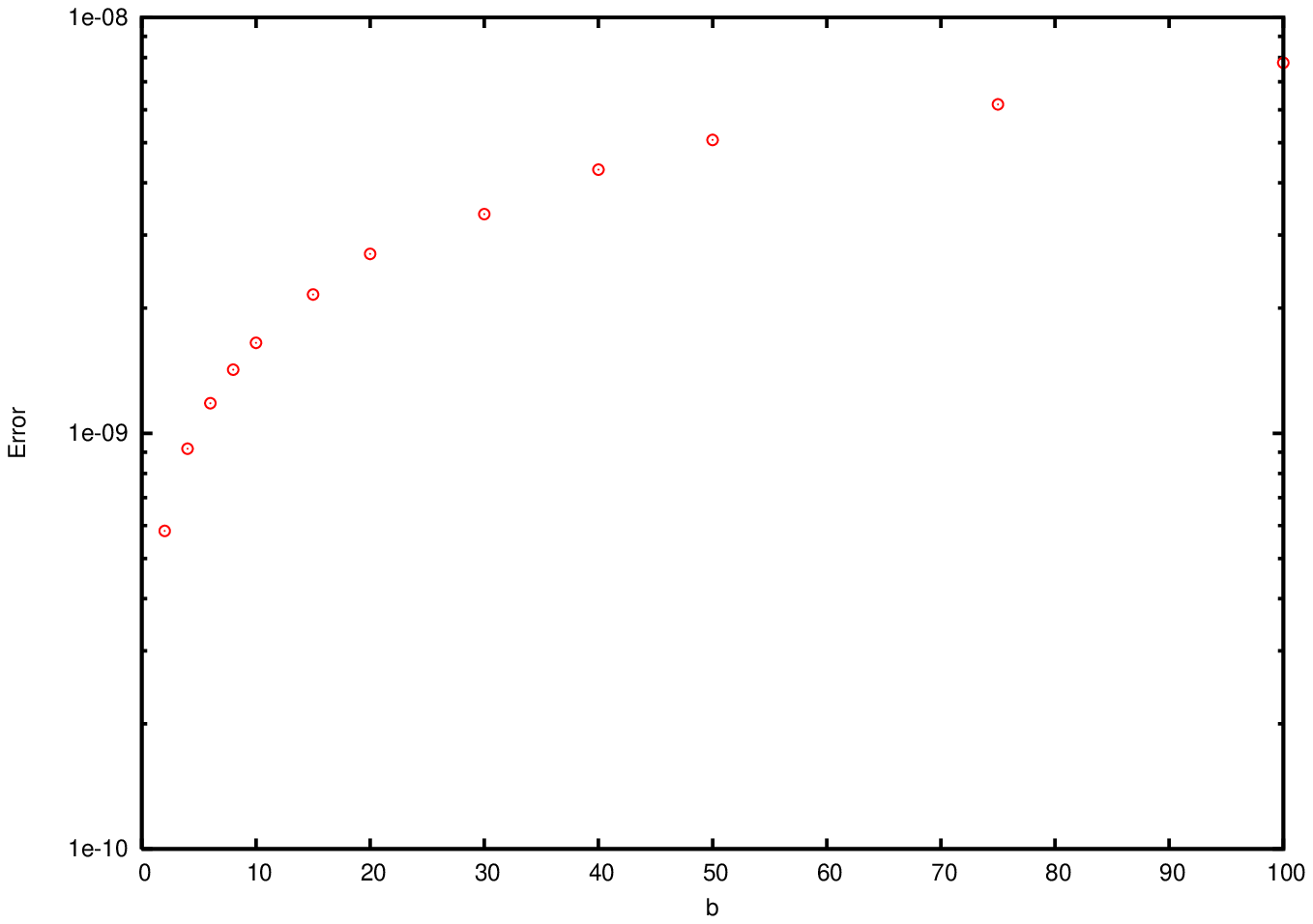}
\caption{ Error as a function of $\blength$, the number of conformal 
maps in a block. The data shown uses $N=100,000$, 
$\nterms=12$, and $L=4$.
}
\label{err_blength}
\end{figure}

Finally we consider the effect of the block size $\blength$. It is not
clear a priori how the error and time will depend on $\blength$. 
Increasing $\blength$ reduces the number of compositions to compute in 
\reff{blockcompose} but it will also increase the radii of convergence 
$R_j$ which will result in the hat power series approximation 
being used less frequently. 
Figure \ref{err_blength} shows that the error does not depend strongly 
on $\blength$ and so $\blength$ should just be chosen to minimize the time. 
This choice depends significantly on $N$. 
(In this figure we have taken $N=100,000$, $\nterms=12$, and $L=4$.)

\begin{figure}[htb]
\includegraphics{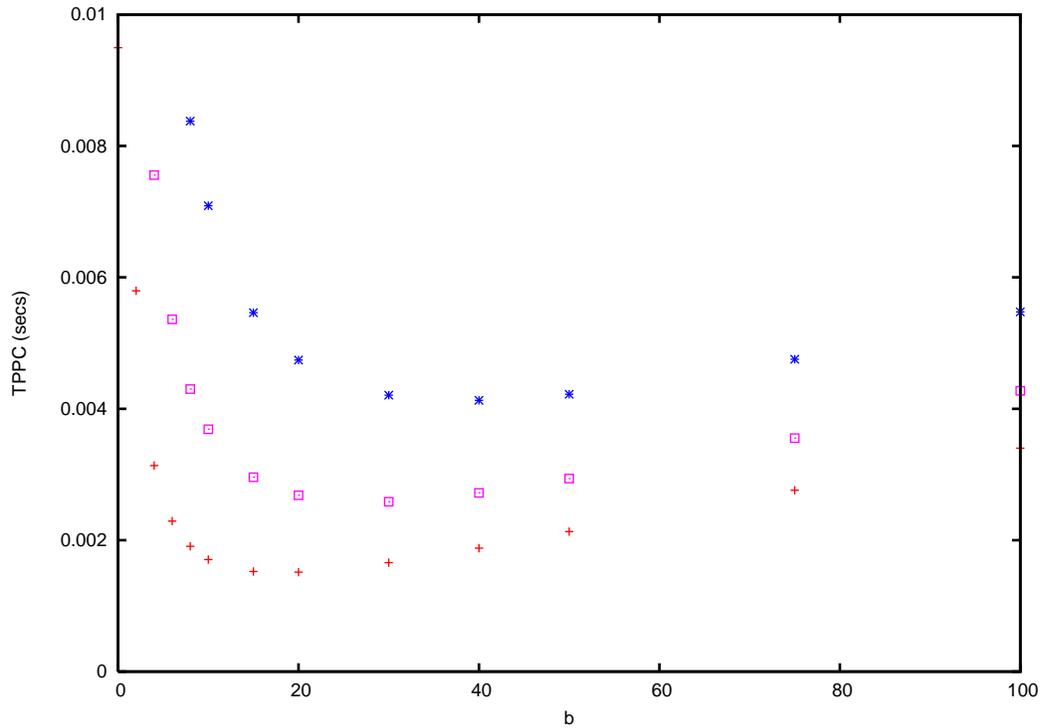}
\caption{ Time per point computed as a function of $\blength$. The three
curves shown from bottom to top are $N=20,000$, $50,000$, and $100,000$. }
\label{time_blength}
\end{figure}

\begin{figure}[htb]
\includegraphics{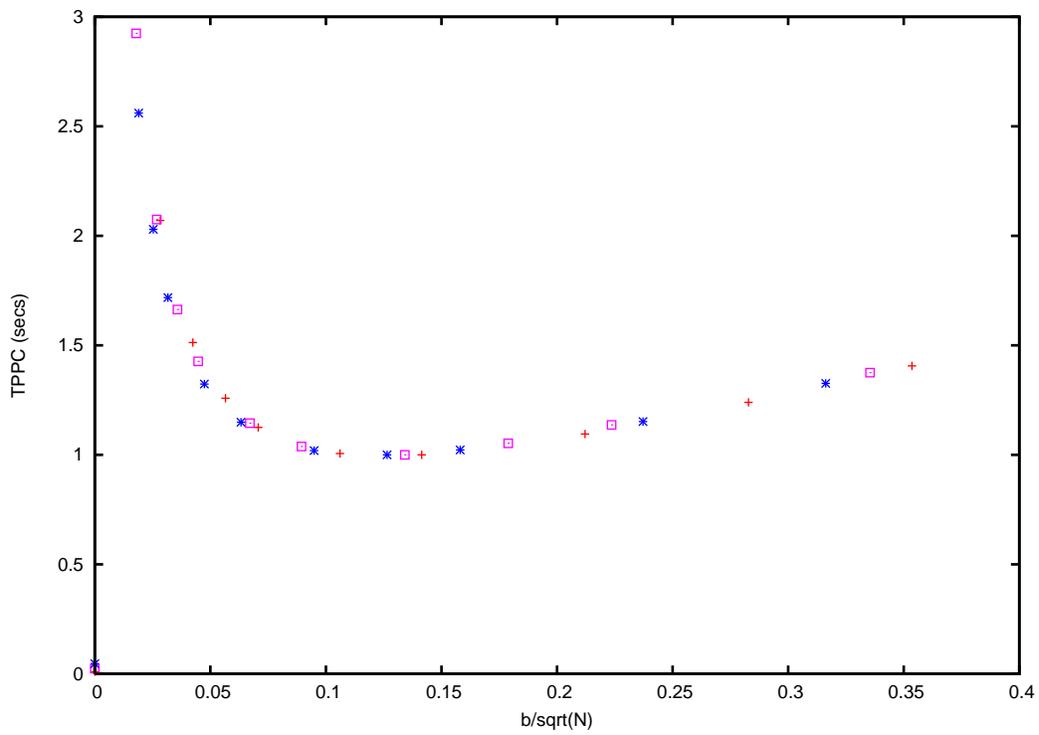}
\caption{ Scaling plot for time per point computed as a function of $\blength$}
\label{time_blength_scaling}
\end{figure}

Figure \ref{time_blength} shows the time as a function of $\blength$ for three
values of $N$. (For all three curves $n=12$ and $L=4$.)
The curves have similar shapes, suggesting some sort of 
scaling. In figure \ref{time_blength_scaling} we plot the same data 
but now divide the time by the minimum time for that value of $N$ and 
divide $\blength$ by $\sqrt{N}$. The resulting three curves 
collapse nicely. This indicates that the optimal value of $\blength$
is roughly proportional to $\sqrt{N}$. 
We have found that the optimal value is well approximated by 
$\blength=0.12 \sqrt{N}$.

We have found that the error is typically well below $L^{-\nterms}$.
To study which value of $\nterms$ is optimal we do the following. 
We fix a value of $\nterms$ and then choose $L$ so that 
$L^{-\nterms}=10^{-6}$. (The choice of $10^{-6}$ is ad hoc.)
We then study the time per point as a function of $N$. Figure 
\ref{time_blength_scaling} shows the resulting plots for 
$\nterms=8, 10, 12, 14$. For the large values of $N$ there is little 
difference between $\nterms=10,12,14$, but they are significantly better
than $\nterms=8$. 

\begin{figure}[htb]
\includegraphics{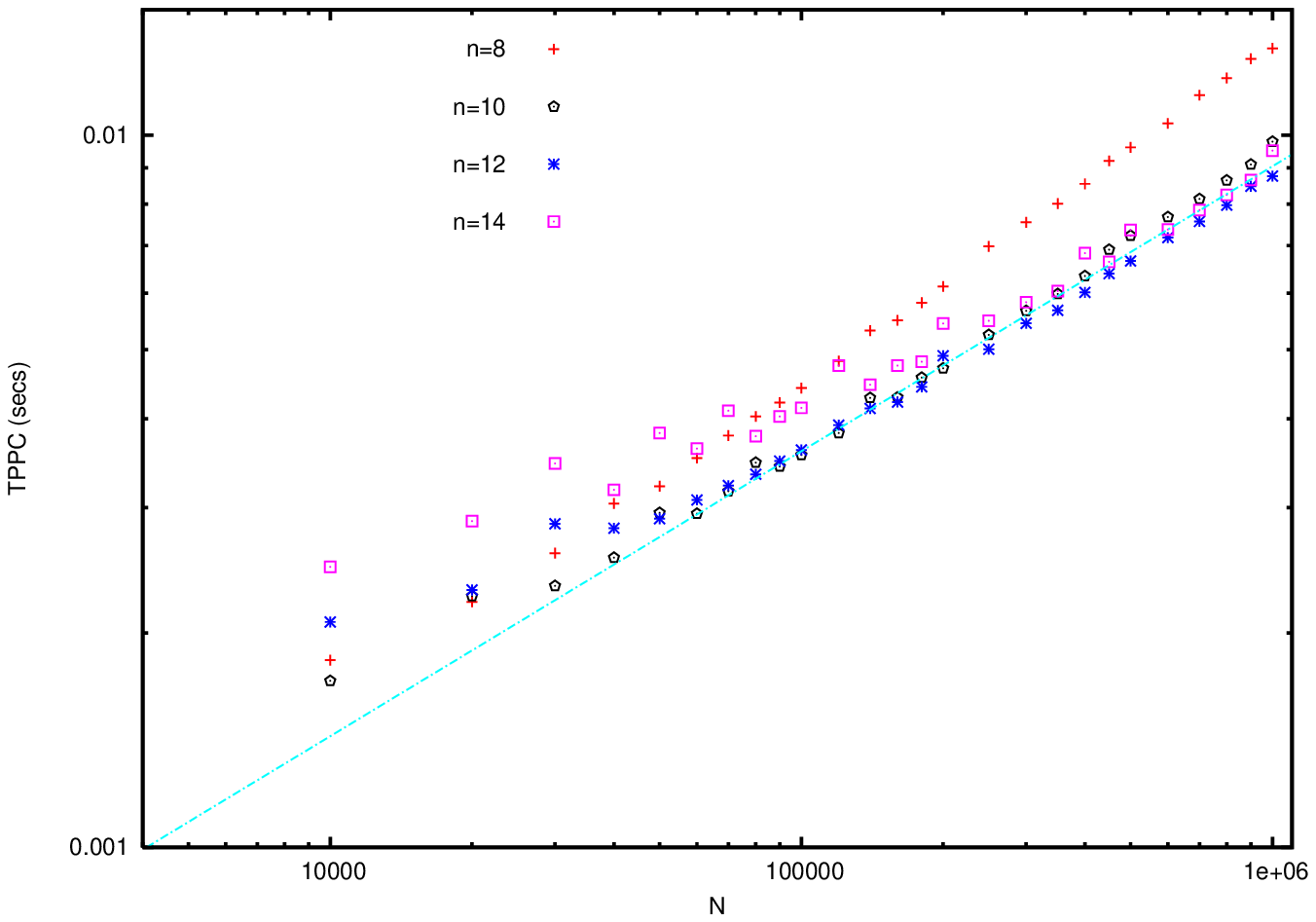}
\caption{ Time per point computed as a function of $N$.
The four curves correspond to $\nterms=8,10,12,14$. The line shown has 
slope $0.4$. 
}
\label{tppc}
\end{figure}

We have carried out the same simulations and generated the same 
plots for $\kappa=6$. Qualitatively the curves are the same. 
The difference between the two values of $\kappa$ 
varies with the choices of the three parameters, but to a very 
crude approximation we have found that for $\kappa=6$ the algorithm 
is about $20\%$ slower, and the error is about twice as large.
(The error for $\kappa=6$ is still usually less than $L^{-n}$.)
The optimal value of $\blength$ for $\kappa=6$ is a bit smaller.
It is better approximated by $\blength=0.1 \sqrt{N}$.
The analog of figure \ref{tppc} for $\kappa=6$ is virtually 
indistinguishable from the figure shown in which $\kappa=8/3$. 
In particular, the time per point is approximately
$O(N^{0.4})$. Taking $\nterms$ to be $10,12$ or $14$ give 
similar results, all significantly better than $\nterms=8$.  

%

\begin{appendix}
\section{Particular hat power series}

In this appendix we give the hat power series needed for the 
two particular discretizations we have discussed. 
For the approximation that uses slits at an angle $\alpha \pi$, we
need to compute the power series of $\hat{h}$ where 
$h(z)= (z+x_l)^{1-\alpha} (z-x_r)^{\alpha}$. 
We have 
\be 
\hat{h}(z)= z \, (1+x_l z)^{-(1-\alpha)} \, (1-x_r z)^{-\alpha} 
\ee
The power series of the last two factors are given by the formula
\be
(1-cz)^{-\alpha} = \sum_{k=0}^\infty \, 
{\alpha (\alpha+1) \cdots (\alpha+k-1) \over k!} \, c^k \, z^k
\label{frac_ps} 
\ee

For the approximation that uses vertical slits, we need to compute
the hat power series of $h(z)=\sqrt{z^2-4 t} +x$.
First consider $g(z)=\sqrt{z^2-4 t}$. We have
\be
\hat{g}(z)={z \over \sqrt{1-4 t z^2}}
= z \sum_{k=0}^\infty \, {1 \cdot 3 \cdot 5 \cdots (2k-1) \over k!} 
2^k t^k z^{2k}
\ee
where the power series may be obtained from \reff{frac_ps} with $\alpha=1/2$. 
Noting that $h=f \compose g$ with $f(z)=z+c$, the hat power
series of $h$ is just the composition of the hat power series for 
$f$ and $g$. The series for $f$ is just 
\be
\hat{f}(z) = { 1 \over 1/z+c} 
= { z \over 1+c z}  = z \sum_{m=0}^\infty (-1)^m \, c^m z^m
\ee

\end{appendix}

\bigskip
\bigskip

\noindent {\bf Acknowledgments:}
The Banff International Research Station made possible many fruitful 
interactions. In particular, I learned much of the material in 
section \ref{sectdiscretize} from 
conversations with Steffen Rohde and Don Marshall.
This work was supported by the National Science Foundation (DMS-0201566
and DMS-0501168.)

\end{document}